\documentclass[12pt]{article}
\usepackage{geometry}
\usepackage{color}
\usepackage{amssymb}
\newtheorem{Theorem} {Theorem}
\newtheorem{Proposition}{Proposition}

\title{
{
\normalsize\tt\hfill\jobname.tex}\\
On convergence rate for Erlang--Sevastyanov type models with infinitely many servers\\
{\small In memory and to the 90th anniversary of A.D. Solovyev (06.09.1927--06.04.2001)}}
\author{A.Yu. Veretennikov\footnote{ University of Leeds, UK; National Research University Higher School of Economics, and Institute for Information Transmission Problems, Moscow, Russian Federation, email: a.veretennikov @ leeds.ac.uk. 
This study has been funded by the Russian Academic Excellence Project '5-100' and by the RFBR grant \mbox{17-01-00633$\_$a}. 
}
}

\begin{document}
\maketitle

\begin{abstract}
Polynomial convergence rate to stationarity is shown for extended Erlang -- Sevastyanov's model with variable intensities of service and arrivals.
\end{abstract}

\section{Introduction}

Consider a service system or
a process with
countably many servers and one
incoming flow of ``events'' or customers ``of the same type'' such that 
the incoming flow has {\em intensity} $\lambda(X)$,
which {\em depends} on the number $n$ of the customers in the system,
and some variables $(x^0, x^1, \dots, x^n)\in R^{n+1}$; we will use notations 
\(X=(n, x^0, x^1, \dots, x^n)\) and \(x=(x^0, x^1, \dots, x^n)\),  where \(n\in Z_+\).
Then, a bit non-rigorously, existence of  intensity \(\lambda(X_t)\) means that
\begin{eqnarray*}
P(\mbox{one new customer arrives on $(t, t+\Delta)$} \mid X_t)
= \lambda(X_t) \Delta + o(\Delta),
\end{eqnarray*}
as $\Delta \downarrow 0$; actually, for {\em continuous} intensities this definition is strict but continuity will not be assumed. Here $X_t = (n_t; x^0_t; x^1_t, \ldots, x^{n_t}_t)$ where  \(n_t\) signifies the number of customers in the system at time \(t\), each $x^i_t, \, 1\le i \le n_t$, stands for the elapsed service time of the corresponding customer at time $t$, while $x^0_t$ signifies the time from the most recent {\em arrival.} The latter makes sense mainly for states where \(n_t=0\); however, for simplicity of presentation we keep the same notations for all states; just for \(n\ge 1\), \(x^0\) must coincide with one of the other \(x^i\) values (namely, with the minimal one). Every customer with the elapsed time of service \(x^i\) is being served by the corresponding server also with some intensity \(h(X^i)\); respectively, the vector  \((x^1_t, \cdots, x^n_t)\) consists of all such elapsed times of all different customers currently present in the system; the meaning of \(x^0_t\) was explained earlier.

In such a state space the process $X_t$ is definitely Markov; however, for discontinuous intensities $\lambda(\cdot)$ a justification of existence of the process is needed, which issue was discussed, e.g., in  \cite{Davis} and \cite{zv-mprf}. In particular,  in \cite{Davis} it has been  proved that the process defined in such a way is well-defined and is strong Markov, which will be useful in the sequel. 

The problem addressed in this paper is convergence rate to the stationary regime under appropriate conditions. Recently such convergence rate was studied for the case where \(\lambda\) may only depend on \(n\) -- the number of customers on service -- but not on other continuous variables.  
Here we address a more general case.
The price for this generality is that the explicit formulae for the stationary distribution due to Fortet and Sevastyanov  are not valid any more (but, of course, modelling is available).

Recall for completeness that Erlang's formulae for the stationary regime in the case of exponentially distributed service time (the system \(M/M/\infty\)) are known since the article \cite{Erlang} (``Erlang formulae {\bf B}'' in the sequel literature),
\begin{equation}\label{Er1}
p_n = p_0 \frac{\lambda^n}{\mu^n}, \quad p_0 = 1-\frac{\lambda}{\mu}.
\end{equation}
where \(\mu^{-1}\) is expectation of the service time,
under the assumption of the convergence of the series (\(\sum_n p_n\)), i.e., for  \(\lambda<\mu\). This was extended under appropriate assumptions (see, e.g., \cite[Ch.4, \S 4-5]{Karlin-Taylor?}) to the case where \(\lambda= \lambda_n\) may depend on the current number of customers \(n\) currently in the system,
\begin{equation}\label{Er2}
p_n = P(\mbox{\(n\) customers in the system})= p_0 \frac{\prod_{i=0}^{n-1}\lambda_i}{\mu^n}, 
\end{equation}
where 
\[
\quad p_0 = \left(\sum_{n=0}^{\infty}
\frac{\prod_{i=0}^{n-1}\lambda_i}{\mu^n}\right)^{-1},
\]
and 
\(\prod_{i=1}^{0}\lambda_i \equiv 1\), and to a  more general situation where {\em intensity} \(\mu\) may also depend on \(n\) (in which case \(\mu^n\) in (\ref{Er2}) should be replaced by \(\prod_{i=1}^{n}\mu_i\)). 
Similar formulae also hold true for finitely many servers in the system under the condition that the customers arriving while all servers are busy become lost.

Fortet \cite{Fortet} for a finite number of servers case  gave the density of a stationary distribution for a general -- {\em non-exponential} -- service time distribution possessing a density \(g(x) = G'(x)\) where \(G\) is the distribution function of service time with a finite mean value \(\displaystyle \mu^{-1}:= \int_0^\infty t\,dG(t) \equiv  \int_0^\infty (1-G(t))\,dt\):
\begin{equation}\label{Fo3}
p(n; x) = p(n; x^1, \dots, x^n) = p_0\,  \prod_{i=0}^{n-1} \lambda_i (1-G(x^i)), \quad n\ge 1, 
\end{equation}
where \(p_0\) is the normalizing constant; in the case of infinite number of servers it is given by the same formula as in (\ref{Er2}).
Sevastyanov \cite{Sevastyanov} extended and strengthened this result to the case without assumption on the existence of a density \(g\)  and for the first time has proved {\em convergence in total variation} of the non-stationary system to its stationary regime using his version of the ergodic theorem for Markov processes with general state spaces. This was a breakthrough and by this reason the author's view is that the name ``Erlang--Sevastyanov systems'' is appropriate here. In \cite{Sevastyanov}  
the number of servers was assumed finite and waiting in a queue unavailable, as in \cite{Fortet}. Further, there were several extensions of this result 
to the case of infinitely many servers \cite{Matthes2, Matthes, Schassberger, Ver-mn}. In all these cases starting from \cite{Fortet}, in particular situations with a constant intensity \(\mu\)  the explicit  formula (\ref{Er2}) holds true, although, the process which equals the number of customers in the system at time $t$ is not Markov. 

The paper \cite{Ver-mn} was based on the MSc project by the author on the topic suggested by Professor A.D.~Solovyev.  The ultimate goal of that time -- to estimate rates of convergence for Erlang--Sevastyanov systems -- was partially realised only recently \cite{Ver-QS} where a polynomial convergence for such systems  with finitely or infinitely many servers has been established under certain assumptions on the intensities; 
a bit earlier in \cite{Ver-Kelb} exponential rate was proved under more severe  than in
\cite{Ver-QS} assumptions. The problem of evaluating convergence rates without assuming anything about intensities apparently remains open. At the same time, neither in \cite{Sevastyanov}, nor in \cite{Ver-mn} existence of intensities was assumed. 
Note that there are also general results about exponential or polynomial convergence for Markov processes and for {\em regeneration} processes to stationarity under the assumption of appropriate exponential or polynomial  recurrence along with a ``local mixing'' or regeneration \cite{Asmussen, Borovkov, Kalashnikov, Kalashnikov2, menshikov95, Thorisson}, et al. So far, none of these general results cover directly  the particular  Erlang--Sevastyanov type systems and their convergence rates and especially the setting with a ``more general'' dependence of \(\lambda(\cdot)\) of the ``whole state'' of the process.   In any case, the results proposed below do not follow from earlier results, including those from \cite{Ver-QS}.  On the other hand, close results for single--server systems were recently studied in 
\cite{Ver-ait, Ver_arxiv2016, zv-mprf},  
and some ideas from the latter articles will be used in the sequel. 
A general dependence of intensity of service \(h(\cdot)\) on all coordinates of the process (assuming that the description of this process makes it Markov) is allowed. 
Initially, the result of this paper was announced in \cite{Ver_arxiv2015}. The literature in the references is by all means not complete; more references can be found in \cite{Ver-QS}.

The paper consists of Introduction, Main  result, Proof of main result.

\section{Erlang -- Sevastyanov type system: main result}
Let us introduce the state space ${\cal X}$: it is a union of  countably many subsets,
$$
{\cal X}= \bigcup_{n=0}^\infty \{(n,x^0, \ldots, x^n), \; x^0, \ldots x^n \ge 0\}; 
$$
denote also ${\cal X}_n := \{(n,x^0, \ldots, x^n), \; x^0, \ldots x^n \ge 0\}$.
To any $n\ge 0$ there correspond $n+1$ non--negative 
coordinates $(x^1, \ldots, x^n)$ and $x^0$, which signify, respectively, the elapsed times
of service of all existing $n$ customers ($x^1, \ldots, x^n$) and the time from the last arrival ($x^0$), including for $n=0$.
It is convenient to assume that a newly arrived customer
with number $n+1$ gets a coordinate $x^i=0$ for any $i=1, \ldots,
n+1$ with equal probabilities. For state \(X=(n, x^0, \ldots, x^n)\), denote \(x = (x^0, \ldots, x^n)\), and  \(n(X) = n(x) = n\) (i.e., we do not distinguish \(n(X)\) and \(n(x)\)).
To establish convergence rate, we will use a characteristic
of the service time distribution called intensity of service,
$$
h(t) : = \frac{g(t)}{1-G(t)}, \quad t \ge 0, \quad
g(t)  = G'(t).
$$
If the intensity function is constant, it means an exponential
distribution of the service time.
With a convention \(\sum_{j=1}^{0}  = 0\), denote for \(X=(n,x)\),
\[
V_{m,a}(X) :=
\left(\sum_{j=1}^{n(X)}(1+x^j)^{m}\right)^{a}, \quad \Lambda:=\sup_{X \in {\cal X}} (\lambda(X)/(n(X)\vee 1)),
\]
\[
\underline{\lambda}_0:= \inf_{Y \in {\cal X}_0} \lambda(Y), \quad 
\bar \lambda_n := \sup_{X \in {\cal X}_n}\lambda(X), \quad n\ge 0.
\]
Note that \(V_{m,a}\) as a function does not depend on \(x^0\) and that this function is symmetric with respect to $(x^1, \ldots, x^n)$. Notation \(a\vee b\) stands for \(\max(a,b)\). The process \(X_t\) is assumed right-continuous.

\begin{Theorem}\label{thm1}
Let there exist $C_0, D>0$, 
$m>1$, $a >1$ and \(\ell > 0\) 
such that 
\begin{equation}\label{ash}
 \frac{C_0}{1+t} \le h(t) \le D, \quad t\ge 0,
\end{equation}
\begin{equation}\label{ash2}
0 < \underline{\lambda}_0 \le \Lambda < \infty,
\end{equation}
and 
\begin{equation}\label{aslambda}
C_0 > \left(a+ \frac{\ell+1}{m}\right)\left( m + \Lambda 2^{a+ \frac{\ell+1}{m}}\right).
\end{equation}
Then for any \(0<k<\ell\), if $C_0$ is large enough, then   
there exist constants $C, C_1>0$  such that for every $X_0 \in {\cal X}$ and \(t\ge 0\),
\begin{equation}\label{esti}
\|\mu^{X_0}_t - \mu \|_{TV} \le
\frac{C(V_{m,a+\ell/m}(X_0)\vee 1+C_1)}{(1+t)^{k+1}}, 
\end{equation}
where \(\mu^{x}_t\) is the distribution of \(X_t\) with the initial value \(x\), \(\mu\) is the unique stationary measure of the process, and \(\|\cdot \|_{TV}\) is the total variation distance.  
\end{Theorem}
{\bf Remark 1}. 
The condition (\ref{aslambda}) used here is a bit more precise than in \cite{Ver-QS}; probably it may allow some further  generalizations. For any  \(k>0\), the bound (\ref{esti}) with {\em some} \(m\) and \(a\) is valid for  any  \(C_0\) which is large enough. 
Recall that strong Markov property holds true for our system, see  \cite{Davis}. 
Uniqueness of stationary distribution as well as its existence is a part of the statement of the Theorem.  

~

\noindent
{\bf Remark 2}.
The constant $C, C_1$ in (\ref{esti}) are not just finite but do admit some effective estimate, which will be seen in the proof. In the earlier versions including \cite{Ver-QS}  similar constants in the main estimate -- e.g., in the main inequality in \cite{Ver-QS} -- also admit certain effective estimates in principle. 
Nonetheless, in the present paper we show step by step how such estimates may be achieved. There is a hope that better and more precise bounds may be obtained in the future. 
More that that, there is a hypothesis that the condition (\ref{ash2}) may be relaxed to a condition similar to (\ref{ash}),  
\begin{equation}\label{ash3}
  \frac{C'_0}{1+t} \le\lambda(0,t) \le \Lambda < \infty, \quad t\ge 0, 
\end{equation}
with some $C'_0$. Yet, this seems more technical and we do not pursue this goal here. 

\section{Proof of Theorem \ref{thm1}}

\newcounter{step}

{\bf 0}. \label{step0}
\refstepcounter{step}
We will be using notation $*$ for an arbitrary value of any coordinate of the process: e.g., $(1,*, 0) = (X = (1,x^0, 0)\in {\cal X}: \, n(X)=1, x^0\ge 0)$. The starting idea is to construct a Lyapunov function and to apply coupling method. The coupling algorithm to be used here will be a ``two-step'' one, which makes it, perhaps, just a little bit unusual. Eventually, we are going to show that roughly speaking the state $(1,0,0)$  may be considered as a ``generalised regeneration'', with some uniformly bounded from above polynomial moments for the distributions of the length between the generalised regeneration periods. (The word ``generalised'' here means that it is necessary to change probability space to make it work.) Due to the strong Markov property, this will suffice for the proof. 

~

Denote 
\[
\tau_{0}:= \inf(t\ge 0: \; X_t \in {\cal X}_0 \equiv \{(0,*)\}), \quad 
\tau_{01}:= \inf(t > \tau_0: \; X_t = (1,0,0)), 
\] 
and 
\[
\tau_{010}:= \inf(t > \tau_{01}: \; X_t = (0,*)).
\]

~

The {\em two-step coupling} idea is to consider two versions of the process, one of them stationary (existence of which is yet to be established), and to couple {\em firstly the first} components of the two processes waiting when they both are simultaneously in the set ${\cal X}_0$, and secondly coupling their remaining components at their transition  from $(0,*)$ to $(1,0,0)$.  

~

\label{step1}
\refstepcounter{step}
\noindent
{\bf \ref{step1}.}  
To estimate $\tau_0$, we construct a Lyapunov function only for the variables $(n, x^1, \ldots, x^n)$, i.e., ignoring the component $x^0$. We want to show that the process $X_t=(n(X_t), X^0_t, X^1_t, \ldots, X^n_t)$ with probability one hits the set \(\{X:\, n(X) = 0\} \equiv {\cal X}_0\), being also positive recurrent uniformly in the variable $x^0$. (Note that given the assumptions this is reasonable, since  only \(\lambda_0(\cdot)\) depends on \(x^0\) essentially, and \(\lambda_0(\cdot)\) itself is bounded away from zero and from infinity.) So, informally  speaking, we may arrange {\em coupling} of two versions of the process -- the original one ($X$) and the stationary one  (say, $Y$) -- on their joint jump from  $n=0$ to $n=1$ (or, more precisely, from the set  ${\cal X}_0 \times {\cal X}_0$ to state $(1,0,0,1,0,0) \in {\cal X}_1 \times {\cal X}_1$). Note, however, that existence of a stationary measure itself is yet to be shown; this will be addressed at step 
\ref{step_stat_measure_i} of this proof which consists of \ref{step10} major steps (this preliminary zero one is not counted): the first part is  fairly close to the calculus from \cite{Ver-QS}, while the second one is  a ``simple coupling'' on the passage of the system ``from \(n=0\) to \(n=1\)''. 
For the convenience of the reader and by the suggestion of the referee we keep most of the calculus trying to minimize references which would make it necessary to read another paper (\cite{Ver-QS}), even though it increases the volume to some extent.  
The constants \(C, C'\), etc. in the calculus may change from line to line. 

For $X=(n,x^0,x^1,\ldots, x^n)$ with $n\ge 1$ and for $1\le j\le n$ denote
\[
X^{j,-}:=(n-1,x^0,x^1, \ldots, x^{j-1}, x^{j+1}, \ldots, x^n)
\]
(i.e., the component \(x^j\) has been dropped, which automatically decreases the value of \(n\) by one unit), and for $1\le j\le n+1$,
\[
X^{j,+}:=(n+1,x^0,x^1, \ldots, x^{j-1},0, x^{j}, \ldots, x^n)
\]
(i.e., a new component \(0\) has been added between the old \(x^{j-1}\) and \(x^{j}\), which automatically increases the value of \(n\) by one unit). For $n=0$ we have $X=(n,x^0)$ with some $x^0\ge 0$, and then (\(j=1\))
\[
X^{1,+}:=(n+1,x^0,0).
\]
The element $X^{j,-}$ is not defined for $n(X)=0$. 

Now, the {\em extended generator} of our {\em Markov process} (see, e.g., \cite{Davis}) reads, 
\begin{eqnarray}\label{gen}
Lg(X) = Lg(n,x) 
= \lambda_{}(X)\frac1{n(x)+1}\sum_{j=1}^{n(x)+1}[g(X^{j,+})-g(X)] 
 \nonumber \\ \\ \nonumber  
+ 1(n(x)>0)\sum_{j=1}^{n(x)}h(x^j)[g(X^{j,-})-g(X)]
+ \sum_{i=0}^{n(x)}\frac{\partial}{\partial x^i}g(n,x).
\end{eqnarray}
It would be a ``honest'' generator of the Markov process -- see \cite{Dynkin} -- under the assumption of continuity of all intensities \(\lambda\) and \(h\). 
A corresponding {\em Dynkin's formula} -- which is, actually, the definition of  extended generator (except that we do not state the exact domain of this operator, but only some sub-domain, which is sufficient for our goals) -- has a form, 
\begin{equation}\label{dfor1}
\mathbb E_{X_0}g(X_t) = g(X_0) 
+ \int\limits_0^t \mathbb E_{X_0} Lg(X_s)\,ds,
\end{equation}
for any function $g$ from the class of bounded $C^1_b$ functions (i.e., with bounded continuous derivatives in any $x^i, \, 0\le i\le n(X)$).
The {\em proof} of this formula for $L$ given in (\ref{gen}) follows from the ``complete expectation'' formula (by analogy with the complete probability formula); some details in a slightly different setting may be found in \cite{zv-mprf}; see also \cite{Davis}.  
Equivalently, the equation (\ref{dfor1}) may be stated as follows: the difference
\begin{equation}\label{mart}
g(X_t) - g(X_0) - \int\limits_{0}^t Lg(X_s)\,ds
\end{equation}
is a martingale (see, e.g., \cite{Davis}). (We note that because of the jumps the accurate writing should have used $Lg(X_{s-})$ under the integral here, but since integration is performed with respect to the Lebesgue measure, it is equivalent to  $Lg(X_{s})$.)

Using just bounded functions is often not sufficient. It follows further that for {\em locally bounded}  $C^1$ functions (i.e., with one continuous derivative in all components $(x^0, \ldots, x^n)$) the same expression (\ref{mart}) is a {\em local martingale} (i.e. a martingale stopped by some appropriate ``localizing sequence'' of stopping times, see, e.g., \cite{RLASH}). 
For the sequel, note that another form of Dynkin's formula for bounded functions $f(t,X)$ from the class $C^1_b$ in all $x^i$ and in $t$ reads (\(f_s\) denotes \(\partial f/\partial s\)), 
\[
\mathbb E_{X_0}f(t,X_t) = f(0,X_0) + \int\limits_0^t \mathbb E_{X_0} (f_s(s,X_s) + Lf(s,X_s))\,ds.
\]
Equivalently, it can be stated that for $f$ from  the specified class of functions the difference
\[
f(t,X_t) - f(0,X_0) - \int\limits_{0}^{t} (f_s(s,X_s) + Lf(s,X_s))\,ds
\]
is a martingale. Similarly to the case of functions $g(X)$, it follows that for {\em locally bounded} $C^{1}$ functions $f(t,X)$ with locally bounded first derivatives with respect to $t$ and all $x^i$ the latter difference is a {\em local martingale}.

~

\label{step2}
\refstepcounter{step}
\noindent
{\bf \ref{step2}}.
Let us show that $V_{m,a}$ may serve as a Lyapunov function. 
For $X_t\not\in {\cal X}_0$ and the constants $m$ and $a$ satisfying a weakened version of the standing assumption (\ref{aslambda}),
\begin{equation}\label{less}
C_0 > a(m +\Lambda 2^a), 
\end{equation}
and with $n = n(X_t)$ and with $M_t$ being some  local martingale we have (cf. with \cite[Lemma 2]{Ver-QS}),
\begin{eqnarray}\label{Ito}
& dV_{m,a}(X_t) 
= \lambda_n(X_t) \,
\left(\left(1+\sum_{j=1}^{n}(1+X_t^j)^{m}\right)^{a}
- \left(\sum_{j=1}^{n}(1+X_t^j)^{m}\right)^{a}\right)
dt
 \nonumber\\ \nonumber\\ \nonumber
& - \sum_{i=1}^n h(X^i_t)
\left(-\left(\sum_{1\le j\le n, \, j\not=i}^{}(1+X_t^j)^{m}\right)^{a}
+ \left(\sum_{j=1}^{n}(1+X_t^j)^{m}\right)^{a}
\right)\,dt
 \\ \nonumber\\ \nonumber
&+ \sum_{i=1}^n \left(
\left(\sum_{j=1}^{n}(1+X_t^j + dt)^{m}\right)^{a}
- \left(\sum_{j=1}^{n}(1+X_t^j)^{m}\right)^{a}
\right)
+ dM_t
 \\ \nonumber \\ 
& \equiv (I_1 - I_2 + I_3)dt+ dM_t. \end{eqnarray}
Due to the assumption (\ref{ash}) we get,
\begin{eqnarray*}
I_2 \ge C_0 \sum_{i=1}^n
(1+X^i_t)^{-1}
(1+X^i_t)^m
V_{m,a-1}(X_{t})
= C_0 V_{m-1,1}(X_{t})
V_{m,a-1}(X_{t}). 
\end{eqnarray*}
Hence, we can see that, at least, the ``main term'' \(-I_2\) is negative for \(X_t\not\in {\cal X}_0\).

Our next task is to show that $I_1$ and $I_3$ are dominated by
$I_2$. Then it would imply that the stationary measure
integrates some polynomial. This, in turn, would allow to
extend our Lyapunov function so as to include some
multiplier that depends on time (see step 2 below). The latter would provide for some $k>0$ a crucial bound
$
E_x \tau_0^{k+1} < \infty
$
along with its quantitative version (see the estimate (\ref{kp1}) below). Finally, the similar inequality for \(\tau_{01}\) would imply coupling
between the original process and its stationary version (but not an immediate regeneration at \(\tau_0\), unlike in \cite{Ver-QS}), which
would mean a certain rate of convergence to the stationary regime.

We estimate (cf. with \cite[Lemma 3]{Ver-QS}),
\begin{eqnarray*}
I_1=\lambda_n(X_t) \,
\left(\left(1+\sum_{j=1}^{n}(1+X_t^j)^{m}\right)^{a}
- \left(\sum_{j=1}^{n}(1+X_t^j)^{m}\right)^{a}\right)
\le \bar\lambda_n a 2^a V_{m,a-1}(X_t). \end{eqnarray*}
Due to the inequality $n  \le V_{m-1,1}(X_{t})$ we find,
\begin{eqnarray*}
 I_1 \le
\bar\lambda_n a 2^a V_{m,a-1}(X_t)
\le
\Lambda n \, a   2^a\, V_{m,a-1}(X_{t})
 \le \Lambda \, a  2^a\, V_{m-1,1}(X_{t}) V_{m,a-1}(X_{t}).
\end{eqnarray*}
Further, we have,
\begin{eqnarray*}
 I_3 = \sum_{i=1}^{n} 
\left(a \left(\sum_{j=1}^{n}(1+X_t^j)^{m}\right)^{a-1}
m \sum_{j=1}^{n}(1+X_t^j)^{m-1}\right)
 = a m \, V_{m-1,1}(X_{t}) V_{m,a-1}(X_{t}).
\end{eqnarray*}
Notice that both the estimate for $I_1$ and the expression for $I_3$ are comparable
with the estimate for $I_2$. Overall,
\begin{eqnarray*}
I_1 - I_2 + I_3 \le -
(C_0 - \Lambda a 2^a - m a) \, V_{m-1,1}(X_{t}) V_{m,a-1}(X_{t}) < 0,
\end{eqnarray*}
for $X_t\not\in {\cal X}_0$. 
So, we get,
\begin{eqnarray}\label{fatou1}
{\mathbb E}_X V_{m,a} (X_{t \wedge \tau_0})
+ (C_0 - a(m +\Lambda 2^a)) {\mathbb E}_X \int\limits_0^{t \wedge \tau_0}
V_{m-1,1}(X_{s}) V_{m,a-1}(X_{s}) \, ds
\le V_{m,a} (X),
\end{eqnarray}
and, as $t\to \infty$, under the assumption of  (\ref{less}), by Fatou's lemma  we obtain,
\begin{eqnarray}\label{fatou2}
{\mathbb E}_X V_{m,a} (X_{\tau_0})
+ (C_0 - a(m +\Lambda 2^a)) {\mathbb E}_X \int\limits_0^{\tau_0}
V_{m-1,1}(X_{s}) V_{m,a-1}(X_{s}) \, ds
\le V_{m,a} (X).
\end{eqnarray}
A formal justification of (\ref{fatou1}) uses a localising sequence from the definition of a local martingale. In our case, let us denote \(T_R:= \inf(t\ge 0: \, n(X_t) + \max_{0 \le i \le n(X_t)}x^i_t \ge R)\) for any \(R>0\). Then, all term in the version of Dynkin's formula
\begin{equation}\label{dfo}
{\mathbb E}_X V_{m,a} (X_{t \wedge \tau_0 \wedge T_R})
= {\mathbb E}_X \int\limits_0^{t \wedge \tau_0 \wedge T_R}
LV_{m,a}(X_{s})\, ds
+ V_{m,a} (X),
\end{equation}
are bounded and the identity (\ref{dfo}) itself follows, e.g., from the ``complete expectation'' formula mentioned earlier. 
In other words, the process 
\[
M_{t \wedge \tau_0 \wedge T_R} : = 
V_{m,a} (X_{t \wedge \tau_0 \wedge T_R})
- \int\limits_0^{t \wedge \tau_0 \wedge T_R} LV_{m,a}(X_{s})\, ds
- V_{m,a} (X)
\]
is a stopped martingale. 
So, by the bounds on all terms of \(LV(X_s)\) above, we get 
\begin{eqnarray}\label{fatou1stopped}
{\mathbb E}_X V_{m,a} (X_{t \wedge \tau_0 \wedge T_R})
+ (C_0 - a(m +\Lambda 2^a)) {\mathbb E}_X \int\limits_0^{t \wedge \tau_0 \wedge T_R}
V_{m-1,1}(X_{s}) V_{m,a-1}(X_{s}) \, ds
\le V_{m,a} (X),
\end{eqnarray}
and now (\ref{fatou1}) follows from (\ref{fatou1stopped}) by Fatou's lemma as \(R\to\infty\), as required, since \(T_R \to \infty\) (because the first component of the process \(X\) has jumps \(\pm 1\), other components increase between their exclusively negative jumps with rate \(+1\), and because intensity of arrivals satisfy the condition \(\bar \lambda_n \le \Lambda \,n\), \(n\ge 1\), so that neither the first nor any other component of the process may achieve infinity over a finite time with a positive probability, cf. \cite{Karlin-Taylor?}) and since at any \(t\) the trajectories of the process are continuous almost surely (the latter due to existence of intensities of both arrivals and service). In the sequel in similar places we will drop this standard part related to localising sequence.

~

Further, in particular, it follows that $\mathbb E_X\tau_0 < \infty$ for any $X$.  In the sequel we shall see that also \(\mathbb E_X\tau_{01}<\infty\), which does signify that the process $X$ is positive recurrent (see the details in steps 2--3 below). According to the
Harris--Khasminsky principle -- fully justified by steps 2--3 and \ref{step_stat_measure_i} below  --  
there is an invariant measure $\mu$
(we show in the sequel that it is unique) and $V_{m-1,1}(X) V_{m,a-1}(X)$ is integrable with respect to this measure \(\mu\). Due to the
the elementary inequality (see \cite{Ver-QS})
$$
 V_{m,1}(X)^{(m-1)/m} \le V_{m-1, 1}(X),
$$
(also recall that $V_{m,a} (X_{\tau_0}) = 1$ and
$V_{m,a} (X) V_{m,b} (X) = V_{m,a+b} (X)$), we obtain,
\begin{eqnarray*}
{\mathbb E}_X V_{m,a} (X_{\tau_0})
+ (C_0 - a(m +\Lambda)) {\mathbb E}_X \int\limits_0^{\tau_0}
V_{m,a-1/m}(X_{s}) \, ds
\le V_{m,a} (X).
\end{eqnarray*}
In particular, for any \(t\), 
\begin{eqnarray}\label{inpar}
{\mathbb E}_X V_{m,a} (X_{\tau_0}) \vee 
{\mathbb E}_X V_{m,a} (X_{t\wedge \tau_0}) 
\le V_{m,a} (X).
\end{eqnarray}

~

\label{step3}
\refstepcounter{step}
\noindent
{\bf \ref{step3}}.
Let $V_{m,a,k}(t,X): = (1+t)^k V_{m,a}(X), \; k<\ell.$
Similarly to the above, we have due to Dynkin's formula,
\begin{eqnarray}\label{i1i2i3}
& dV_{m,a,k}(t,X_t)= V_{m,a,k}(t+dt,X_{t+dt}) - V_{m,a,k}(t,X_{t})
  \nonumber \\ \nonumber \\
& = 
(1+t)^{k}\left[I_1 - I_2 +I_3\right]\,dt
+ d\tilde M_t
+ k(1+t)^{k-1} \,  V_{m,a}(X_{t})\, dt
  \nonumber \\  \\ \nonumber 
& \le -(1+t)^k (C_0 - a(m +\Lambda))
V_{m,a-1/m}(X_{t}) \, dt
+ k(1+t)^{k-1} \, V_{m,a}(X_{t})\, dt + d\tilde M_t.
\end{eqnarray}

Now the task is to ensure that the negative part in the
right hand side of the last expression prevails all other terms. We will be using the
inequality established in the step 1 above. The second term may
be split into two parts,
\begin{eqnarray*}
I:= k(1+t)^{k-1} \,  V_{m,a}(X_{t}) \hspace{2cm}
 \\\\
\equiv I\times 1(k(1+t)^{k-1} \, V_{m,a}(X_{t})
\le \epsilon (1+t)^k V_{m,a-1/m}(X_{t}))
 \\\\
+ I\times  1(k(1+t)^{k-1} \, V_{m,a}(X_{t})
> \epsilon (1+t)^k V_{m,a-1/m}(X_{t})).
\end{eqnarray*}
The first part of this term here with '$\le \epsilon$', is dominated by
the main negative expression if $\epsilon < C_0 - a(m + \Lambda 2^a)$.

Let us estimate the second part of this term. We have, for any $\ell>k$,
\begin{eqnarray*}
I\times  1(k(1+t)^{k-1} \, V_{m,a}(X_{t})
> \epsilon (1+t)^k V_{m,a-1/m}(X_{t}))
 \\\\
\le I \times  \frac{(k \, V_{m,a}(X_{t}))^\ell}
{(\epsilon (1+t) V_{m,a-1/m}(X_{t}))^\ell}
= I \times  \frac{k^\ell}{(\epsilon (1+t))^\ell}
V_{m,\ell/m}(X_{t}).
\end{eqnarray*}
Therefore, the second part of the second term does not exceed
\begin{eqnarray*}
k(1+t)^{k-1} \,
\times \frac{k^\ell}{(\epsilon (1+t))^\ell}
V_{m,a+\ell/m}(X_t)^{}.
\end{eqnarray*}
Denote $a': =a+\ell/m$,
and recall that, in any case, \begin{equation}\label{aprim}
C_0 > a'(m + \Lambda 2^{a'}).
\end{equation}
Now, let us collect all terms and their bounds,
integrate and take expectation,
\begin{eqnarray}\label{integrate}
& \displaystyle {\mathbb E}_X V_{m,a,k}(t\wedge \tau_0,X_{t\wedge \tau_0})
+ (C-a(m+\Lambda)-\epsilon) {\mathbb E}_X \int\limits_0^{t\wedge \tau_0}
(1+s)^{k} V_{m,a-1/m}(X_s) \,ds
 \\\\
& \displaystyle  \le V_{m,a}(X)
+ C' \int\limits_{0}^{\infty} {\mathbb E}_X 1(s\le t\wedge\tau_0)
(1+s)^{k-1-\ell} V_{m,a+\ell/m}(X_s)\,ds.
\end{eqnarray}
Here due to (\ref{inpar}) and under the assumption (\ref{aprim}), 
\begin{eqnarray*}
{\mathbb E}_X 1(s\le t\wedge\tau_0)
V_{m,a+\ell/m}(X_s) 
\le V_{m,a+\ell/m}(X). 
\end{eqnarray*}
Thus, 
\[
\int\limits_{0}^{\infty} {\mathbb E}_X 1(s\le t\wedge\tau_0)
(1+s)^{k-1-\ell} V_{m,a+\ell/m}(X_s)\,ds \le C V_{m,a+\ell/m}(X),  
\]
where the constant $C$ clearly admits an effective estimate, as well as $C'$ a few lines earlier. 
Further, due to Fatou's lemma, for 
$k < \ell$ this implies,
\begin{eqnarray*}
{\mathbb E}_X V_{m,a,k}(\tau_0,X_{\tau_0})
+ C' {\mathbb E}_X \int\limits_0^{\tau_0}
(1+s)^{k} V_{m,a-1/m}(X_s) \,ds
\le V_{m,a}(X)
 + C'' V_{m,a+\ell/m}(X).
\end{eqnarray*}
Since $V_{m,a-1/m}(X_s)\ge 1$ for \(s<\tau_0\), we obtain
\[
{\mathbb E}_X \tau_0^{k+1}
\le C V_{m,a}(X)  + C V_{m,a+\ell/m}(X).
\]
or, with one more new constant \(C\), 
\begin{eqnarray}\label{kp1}
{\mathbb E}_X \tau_0^{k+1}
\le C V_{m,a+\ell/m}(X).
\end{eqnarray}
Here also the constant $C$ does allow an effective estimate. Emphasize that so far we have still used a bit relaxed version (\ref{aprim}) (namely, \(C_0 > a'(m + \Lambda 2^{a'})\) with \(a'=a+\ell/m\)) of the standing assumption (\ref{aslambda}).

~

\label{step4}
\refstepcounter{step}
\noindent
{\bf \ref{step4}}.
For the hitting time \(\tau_{01}\) defined earlier 
it follows straight away due to (\ref{kp1}) and from the assumption (\ref{ash2})  that 
\begin{equation}\label{tau01}
{\mathbb E}_{X_{\tau_0}} (\tau_{01} - \tau_0)^{k+1}  \le \int\limits_0^\infty s^{k+1} 
\underline{\lambda}_0 \exp(-\underline{\lambda}_0 s)\,ds 
= \underline{\lambda}_0^{-(k+1)}
\Gamma(k+2) <\infty,
\end{equation}
and, hence, since $V_{m,a+\ell/m}(X)\ge 1$ for $n(X)>0$, 
\begin{eqnarray}\label{kp11}
{\mathbb E}_X \tau_{01}^{k+1}
= 
{\mathbb E}_X {\mathbb E}_{X_{\tau_0}} (\tau_0 + (\tau_{01} - \tau_0))^{k+1} \le C V_{m,a+\ell/m}(X) \quad (n(X)>0).
\end{eqnarray}
Recall for completeness that, according to (\ref{tau01}), for $n(X)=0$ we have ${\mathbb E}_X \tau_{01}^{k+1}\le C=
\underline{\lambda}_0^{-(k+1)}
\Gamma(k+2) <\infty$. Hence, (\ref{kp11}) can be written in a more general form for {\em any} initial value \(X\) as follows, 
\begin{eqnarray}\label{kp110}
{\mathbb E}_X \tau_{01}^{k+1}
= 
{\mathbb E}_X {\mathbb E}_{X_{\tau_0}} (\tau_0 + (\tau_{01} - \tau_0))^{k+1} \le C V_{m,a+\ell/m}(X)\vee 1.
\end{eqnarray}

Moreover, since \(V_{m,a}(X)\equiv 0\) for \(X\in {\cal X}_0\), and so, \(1(\tau_0\le s < \tau_{01})V_{m,a}(X_s)\equiv 0\),  we also get  under the full assumption (\ref{aslambda}) with ``$\ell+1$'', that  
\begin{eqnarray}\label{intl}
& \displaystyle {\mathbb E}_X \int\limits_0^{\tau_{01}}
V_{m,a+\ell/m}(X_{s}) \, ds
= {\mathbb E}_X (\int\limits_0^{\tau_{0}} + \int\limits_{\tau_{0}}^{\tau_{01}})
V_{m,a+\ell/m}(X_{s}) \, ds
 \nonumber \\\nonumber \\
& \displaystyle  = {\mathbb E}_X \int\limits_0^{\tau_{0}} 
V_{m,a+\ell/m}(X_{s}) \, ds
\le C V_{m,a+(\ell+1)/m} (X) \vee 1.
\end{eqnarray}

~

\label{step5}
\refstepcounter{step}
\noindent
{\bf \ref{step5}}.
Now, for the hitting time 
\begin{eqnarray*}
\tau_{010} := \inf(t\ge \tau_{01}: \; X_t = (0,*)), 
\end{eqnarray*}
which will play its role in the sequel, we can use the same estimate (\ref{kp1}) from the previous step, which gives the following bound:  
\begin{eqnarray}\label{kp111}
{\mathbb E}_{X_{\tau_{01}}} (\tau_{010}-\tau_{01})^{k+1}
\le C V_{m,a+\ell/m}((1,0,0)) = C. 
\end{eqnarray}
As an immediate corollary we also get from (\ref{kp111}) and (\ref{kp11}) that 
\begin{eqnarray}\label{kp1111}
{\mathbb E}_X \tau_{010}^{k+1}
\le  C V_{m,a+(\ell+1)/m} (X)\vee 1. 
\end{eqnarray}
Naturally, all constants $C$ here (and earlier) may depend on $k$. 

~

\label{step6}
\refstepcounter{step}
\noindent
{\bf \ref{step6}}.
{\bf Sufficiency of the bound (\ref{kp11}).} For some readers it could be a ``common knowledge'' that the inequality (\ref{kp11})  suffices for the proof of the statement of the Theorem because the state \((1,0,0)\) may be regarded as a regeneration state. 
However, for the author it was not so easy to provide a proper reference which would justify exactly what is stated including the bound (\ref{esti}). In particular, one of the difficulties could be that there is no single regeneration state achievable by the process with a positive probability at any time. As an example, state $(0,0)$ could have potentially served as a regeneration state, but the sojourn time of the process at this state is zero and similarly for all other states. 
So, for the sake of completeness and for the convenience of the reader we offer a version of coupling below. The reader who knows a right reference may skip the remaining paragraphs.

~

\label{step7}
\refstepcounter{step}
\noindent
{\bf \ref{step7}}.
{\bf ``Hitting the set \({\cal X}_0 \times {\cal X}_0\)''}.
Similarly to (\ref{kp1}), it may be justified {\em for the couple of independent  processes}  \((X_t,Y_t)\) with possibly different initial values that
\[
{\mathbb E}_{X_0, Y_0}\bar\tau_{00}^{k+1}<\infty,
\]
with a bound,
\begin{equation}\label{eq00}
{\mathbb E}_{X_0, Y_0}\bar\tau_{00}^{k+1}
\le C\left(V_{m,a+\ell/m}(X_0) + V_{m,a+\ell/m}(Y_0)\right), 
\end{equation}
with any \(k<\ell\), where
$$
\bar\tau_{00} := \inf(t\ge 0: (X_t, Y_t)\in {\cal X}_0 \times {\cal X}_0).
$$
The inequality (\ref{eq00}) can be established similarly to the lines of proving (38--39) in \cite{Ver-QS}. So, let us consider a couple of processes $X_t$ and $Y_t$ on two independent probability spaces where the process $Y_t$ is a stationary version of $X_t$ (with a stationary initial distribution: the existence of this distribution will be {\em independently} justified in step \ref{step_stat_measure_i} in the sequel).  
Consider new Lyapunov functions for the  two component process:
\[
\bar V_{m,a}(X,Y) := V_{m,a}(X) + V_{m,a}(Y),
\]
and 
\[
\bar V_{k,m,a}(t,X,Y) := V_{k,m,a}(t,X) + V_{k,m,a}(t,Y).
\]

It follows from Dynkin's formulae for the couple $(X_t, Y_t)$ and for the function $\bar V_{m,a,k}(t,X_t, Y_t)$ with the restrictions on the parameters given in the assumptions of the Theorem that for $(X_t,Y_t)\not\in {\cal X}_0\times {\cal X}_0$ we have, 
\begin{eqnarray*}
& \displaystyle d \bar V_{m,a,k}(t,X_t, Y_t)
 \\\\
& \displaystyle 
=
(1+t)^{k}\left[I^X_1 - I^X_2 +I^X_3 + I^Y_1 - I^Y_2 +I^Y_3 \right]\,dt
+ k(1+t)^{k-1} \, \bar V_{m,a}(X_{t},Y_t)\, dt
+ d\bar M_t, 
\end{eqnarray*}
with some local martingale $\bar M_t$ and with natural notations 
\begin{eqnarray*}
& \displaystyle I^X_1 =\lambda_n(X_t) \,
\left(\left(1+\sum_{j=1}^{n}(1+X_t^j)^{m}\right)^{a}
- \left(\sum_{j=1}^{n}(1+X_t^j)^{m}\right)^{a}\right), 
 \\ \nonumber\\ \nonumber
& \displaystyle I^X_2 =\sum_{i=1}^n h(X^i_t)
\left(-\left(\sum_{1\le j\le n, \, j\not=i}^{}(1+X_t^j)^{m}\right)^{a}
+ \left(\sum_{j=1}^{n}(1+X_t^j)^{m}\right)^{a}
\right), 
 \\\\
& \displaystyle I^X_3 = 
\sum_{i=1}^{n} 
\left(a \left(\sum_{j=1}^{n}(1+X_t^j)^{m}\right)^{a-1}
m \sum_{j=1}^{n}(1+X_t^j)^{m-1}\right)
\end{eqnarray*}
(former $I_1, I_2$, and $I_3$), and 
\begin{eqnarray*}
& \displaystyle I^Y_1 =\lambda_n(Y_t) \,
\left(\left(1+\sum_{j=1}^{n}(1+Y_t^j)^{m}\right)^{a}
- \left(\sum_{j=1}^{n}(1+Y_t^j)^{m}\right)^{a}\right), 
 \\ \nonumber\\ \nonumber
& \displaystyle I^Y_2 =\sum_{i=1}^n h(Y^i_t)
\left(-\left(\sum_{1\le j\le n, \, j\not=i}^{}(1+Y_t^j)^{m}\right)^{a}
+ \left(\sum_{j=1}^{n}(1+Y_t^j)^{m}\right)^{a}
\right), 
 \\\\
& \displaystyle I^Y_3 = 
\sum_{i=1}^{n} 
\left(a \left(\sum_{j=1}^{n}(1+Y_t^j)^{m}\right)^{a-1}
m \sum_{j=1}^{n}(1+Y_t^j)^{m-1}\right).
\end{eqnarray*}
For any $t$ where both $X_t\not \in {\cal X}_0$ and $Y_t\not \in {\cal X}_0$, the terms $I_2^{X}$ and $I_2^Y$ dominate their counterparts $I^X_1, I_3^X$ and $I^Y_1, I_3^Y$, respectively. Let us inspect the case where one of the processes belongs to ${\cal X}_0$, say, $n(Y_t)=0$. In this case $X_t\not \in {\cal X}_0$, but the term $I^Y_2 = 0$, and we have, 
\[
I^Y_1 = \lambda_0(Y_t) \le \bar\lambda_0, \quad I^Y_3 = 0
\]
(the latter because \(\sum_{i=1}^0 = 0\)). Hence, $n(Y_t)=0$ implies
\[
I^Y_1-I^Y_2+I^Y_3 = \lambda_0(Y_t) \le \bar\lambda_0.
\]
So, the total sum of non-negative terms admits the bound (cf. with (\ref{i1i2i3})), 
\begin{eqnarray*}
I^X_1+I^X_3+I^Y_1+I^Y_3 \le (1+t)^k [(a(m +\Lambda))
V_{m,a-1/m}(X_{t}) + \bar \lambda_0] 
+ k(1+t)^{k-1} \, V_{m,a}(X_{t}), 
\end{eqnarray*}
while the modulus of the only negative term (\(I^X_2\)) equals
\[
C_0 (1+t)^k 
V_{m,a-1/m}(X_{t}).
\]
In the case $n(X_t)=0$ instead, the situation is, of course, symmetric for \(X\) and \(Y\) and the same conclusions hold. 
So, the earlier considerations for one process \(X_t\) in the steps \ref{step2}--\ref{step3} based on the inequality
\[
C_0 > a(m +\Lambda)
\]
remain valid and lead to the desired estimate (\ref{eq00}) if we only check that 
\[
C_0 > a(m +\Lambda) + \bar \lambda_0.
\] 
However, the latter bound does follow from the standing assumption (\ref{aslambda}):
\[
C_0 > \left(a+ \frac{\ell+1}{m}\right)\left( m + \Lambda 2^{a+ \frac{\ell+1}{m}}\right), 
\]
since $a>1$ and, hence,  $C_0> am + 2a\Lambda > 
am + a\Lambda + \bar \lambda_0$.

~

\label{step8}
\refstepcounter{step}
\noindent
{\bf \ref{step8}}.
{\bf Coupling at the ``${\cal X}_0 \mapsto {\cal X}_1$'' passage.} 
Further, from any state \((0,*,0,*) \in {\cal X}_0 \times {\cal X}_0\), coupling is achieved with a positive (bounded away from zero) probability over a unit time interval on the  passage to  state \((1,0,0,1,0,0) \in {\cal X}_1 \times {\cal X}_1\) due to the assumption
\[
0< \underline \lambda_0 
\le 
\bar\lambda_0 <\infty.
\]
Indeed, denote 
\begin{eqnarray*}
\bar \tau_{0011}:= \inf(t\ge 0: \; (X_t,Y_t) = (1,0,0,1,0,0) \; \mbox{after visiting} \; (0,*,0,*)), \\\\
\bar \tau_{001}:= \inf(t\ge 0: \; (X_t,Y_t) \not \in {\cal X}_0\times {\cal X}_0 \; \mbox{after visiting} \; (0,*,0,*)). 
\end{eqnarray*}
It may happen that \(\bar \tau_{001} = \bar \tau_{0011}\), although, in general, \(\bar \tau_{001} \le \bar \tau_{0011}\).   
The moment \(\bar \tau_{0011}\) may be regarded as a moment of coupling, i.e., the moment where the two versions of the process meet up, after which due to the strong Markov property the processes may be considered as equal, of course, after a corresponding change of the probability space. (Clearly, without such a change the moment \(\bar \tau_{001}\) can never occur, i.e., it can equal infinity.) Indeed, after such a change, we may assume that each of the two processes perform jumps up from state \((0, *)\) according to two {\em independent} flows of ``events'', one with intensity \(\underline{\lambda}_0\) and the other due to the ``remainder'' flow with intensity \(\lambda_0(\cdot) - \underline{\lambda}_0\). The point is that the flow with intensity \(\underline{\lambda}_0\) may be regarded as {\em the same} for the two processes and, hence, on this new probability space both processes (more precisely, their equivalents) jump up to state \((1,0,0)\) simultaneously with probability at least
\(p_{01} := (1-\exp(-\underline\lambda_0))\times \exp(-2 (\bar \lambda_0 -\underline\lambda_0))\) on a unit interval of time. 

Note that similarly to the bound (\ref{kp111}) and due to (\ref{eq00}), for the stopping time 
\[
\bar\tau_{00100}:= \inf(t> \bar\tau_{001}: \, (X_t, Y_t) \in {\cal X}_0 \times {\cal X}_0)
\]
the following holds true with some finite constant $C>0$:
\begin{eqnarray}\label{kp1111}
{\mathbb E}_{X_{\bar\tau_{001}}, Y_{\bar\tau_{001}}} (\bar\tau_{00100} - \bar\tau_{001})^{k+1}
\le C  
\sup_{X,Y\in ({\cal X}_0\times (1,0,0))\cup  ((1,0,0)\times {\cal X}_0)}
\bar V_{m,a+\ell/m}(X,Y) = C. 
\end{eqnarray}
Also, similarly to (\ref{tau01}) and due to exactly the same calculus, 
\begin{eqnarray}\label{kp11111}
\sup_{X,Y\in {\cal X}_0} {\mathbb E}_{X,Y} \bar\tau_{001}^{k+1}
\le C = \underline{\lambda}_0^{-(k+1)}
\Gamma(k+2).
\end{eqnarray}

Let us show that after the declared change of probability space, the moment \(\bar\tau_{001}\) satisfies the bound,
\begin{equation}\label{eq10}
{\mathbb E}_{X_0,Y_0} \bar \tau_{001}^{k+1}<\infty,
\end{equation}
and that moreover, for {\bf any} (this is the difference between (\ref{kp11111}) and (\ref{eq10a}) below) \(X_0=(n^1,x), \, Y_0=(n^2,y)\) and for \(k<\ell\) we have, \begin{eqnarray}\label{eq10a}
{\mathbb E}_{X_0,Y_0} \bar \tau_{001}^{k+1} \le 
C \bar V_{m,a+\ell/m}(X_0, Y_0) \vee 1,
\end{eqnarray}
with a new $C>0$. 
The crucial bound (\ref{eq10a}) here is due to (\ref{eq00}) and to the fact that coupling on the passage from ${\cal X}_0 \times {\cal X}_0$ to state $(1,0,0,1,0,0)$ occurs with a positive probability (\(p_{01}\)) over a unit time (see above), with the help of one more geometric like series, as will be shown below. 
Indeed, denote $\bar\tau_{001}=:\bar\tau_{001}(1)$, $\bar\tau_{00100}=:\bar\tau_{00100}(1)$, $\bar\tau_{00100}(0)=:\bar\tau_{00}$, and by induction, 
\begin{eqnarray*}
& \bar\tau_{00100}(i):= 
\inf(t > \bar\tau_{001}(i): \,  (X_t,Y_t) = (0,*,0,*)), 
 \\\\
& \bar\tau_{001}(i+1) := \inf(t > \bar\tau_{00100}(i): \,  (X_t,Y_t) \not \in {\cal X}_0 \times {\cal X}_0), \quad i\ge 1.  
\end{eqnarray*}
At each moment $\bar\tau_{001}(i)$ coupling is possible with a positive probability \(p_{01}\) -- see the definition above -- bounded away from zero. One of these moments will be a ``successful coupling'', i.e., will be  equal to $\bar\tau_{0011}$. (The change of probability space is assumed as prescribed earlier). 
By virtue of the bounds (\ref{eq00}), (\ref{kp1111}) and (\ref{kp11111}) and using the representation for any \(i\ge 1\), 
\[
\bar\tau_{001}(i) = \bar \tau_{00} + (\bar\tau_{001}-\bar\tau_{00}) +  \sum_{j=1}^{i-1} \left((\bar\tau_{001}(j+1) - \bar\tau_{00100}(j))+ (\bar\tau_{00100}(j) - \bar\tau_{001}(j))\right),
\]
we get by induction with some new $C>0$ for any \(i\ge 1\), 
\begin{equation}\label{bdti}
{\mathbb E}_{X,Y}\bar\tau_{001}^{k+1}(i) 
\le C i^k \left(V_{m,a+\ell/m}(X) + V_{m,a+\ell/m}(Y)\right)\vee 1 = C i^k \bar V_{m,a+\ell/m}(X,Y) \vee 1.
\end{equation}
From (\ref{bdti}) we finally get the decisive estimate for the coupling moment \(\bar\tau_{0011}\), 
\begin{eqnarray}\label{L}
& {\mathbb E}_{X,Y}\bar\tau_{0011}^{k+1} 
\le \sum\limits_{i=1}^{\infty}{\mathbb E}_{X,Y}\bar\tau_{001}^{k+1}(i)(1-p_{01})^{i-1}
  \nonumber  \\ \nonumber \\
& \le  C \sum\limits_{i=1}^{\infty} (1-p_{01})^{i-1}
 i^k \bar V_{m,a+\ell/m}(X,Y)\vee 1 \le C \bar V_{m,a+\ell/m}(X,Y)\vee 1.
\end{eqnarray} 

~

\label{step_stat_measure_i}
\refstepcounter{step}
\noindent
{\bf \ref{step_stat_measure_i}}.
{\bf Existence} of some invariant measure \(\mu\) follows from the fact that \((1,0,0)\) is a regeneration state and from (\ref{kp1}) and (\ref{kp11}) (suffices for one component), by the Harris--Khasminsky formula
\begin{equation}\label{reg}
\mu(A) = c \,{\mathbb E}_{(1,0,0)} \int\limits_0^{\tau_{01}} 1(X_t\in A)\, dt)
\end{equation}
(recall that the stopping time \(\tau_{01}\) was defined earlier in the step 0). 
Here  \(c\)  is the normalising  constant. 
For the sequel notice that  by virtue of (\ref{reg})  and (\ref{intl}), 
\begin{equation}\label{hhp}
\int V_{m,a+\ell/m}(x)\, \mu(dx) < \infty, 
\end{equation}
however, for our aims the value of this integral also has to be estimated. 
From (\ref{reg}) (and with the constant $c$ from (\ref{reg})), we have by virtue of (\ref{intl}),
\begin{eqnarray*}
& \displaystyle \int V_{m,\,a+\ell/m}(X)\, \mu(dX) =
c \, {\mathbb E} _{(1,\,0,\,0)} \int\limits_0^{\tau_{01}}
V_{m,\,a+\ell/m}(X_{s}) \, ds
 \\\\
& \displaystyle \le c\,\frac{1}{C_0 - a'(m +\Lambda)} \,V_{m,\,a+(\ell+1)/m} (1,\,0,\,0) = \frac{c}{C_0 - a'(m +\Lambda)}, 
\end{eqnarray*}
with \(a'=a+\ell/m\), because $V_{m,\,a+(\ell+1)/m} (1,\,0,\,0) = 1$.
So, in order to obtain an upper bound for the integral in (\ref{hhp}), it remains to estimate the constant $c$ from above. From (\ref{reg}) with $A=\Omega$ we get, 
\begin{eqnarray*}
\displaystyle c = \frac{1}{{\mathbb E} _{(1,\,0,\,0)}\tau_{01}} 
\le 
\frac{1}{{\mathbb E} _{(1,\,0,\,0)}\tau_{0}+ {\mathbb E} _{(0,\,0)}\tau_{01}} 
\le  \frac{1}{\displaystyle \frac{1}{\bar\lambda_0}+\frac{1}{D}}\equiv \frac{\bar\lambda_0 D}{\bar\lambda_0+D}.
\end{eqnarray*}
This is firstly because the value of $\tau_{0}$ is stochastically minorated by a similar hitting time for a pure death process with a constant intensity $D$ (cf. with the assumption (\ref{ash})) which mean value equals $D^{-1}$, so that  
\[
{\mathbb E} _{(1,\,0,\,0)}\tau_{0} \ge D^{-1}, 
\]
and secondly, because  
the value of $\tau_{01}$ is stochastically minorated by a similar hitting time for a pure birth process with a constant intensity $\bar \lambda_0$ which mean value equals $\bar\lambda_0^{-1}$, so that  
\[
{\mathbb E} _{(0,\,0)}\tau_{01} \ge \bar\lambda_0^{-1}. 
\]

Therefore, we obtain the following bound for the value of the integral in (\ref{hhp}):
\begin{eqnarray}\label{invconst}
 \int\limits_{\mathcal{S}} V_{m,\,a+\ell/m}(X)\, \mu(d X) =
c\, {\mathbb E} _{(1,\,0,\,0)} \int\limits_0^{\tau_{01}}
V_{m,\,a+\ell/m}(X_{s}) \, d s \le 
 \nonumber \\\\ \nonumber
\le 
\frac{\bar\lambda_0 D}{(\bar\lambda_0+D)(C_0 - (a+\ell/m)(m +\Lambda))} =:C_1.\end{eqnarray}
Note that here again the assumption (\ref{aslambda}) was essentially used while deriving (\ref{invconst}), precisely with ``$\ell+1$'', since (\ref{intl}) was applied  for which the condition (\ref{aslambda}) is required.

~

Returning now to (\ref{eq10a}), we can see that it follows by integration that for the distribution \(\mu\) (and with a new constant \(C\)),
\begin{equation}\label{eq10bb}
{\mathbb E}_{X_0, \mu} \bar\tau_{0011}^{k+1}
\le C V_{m,a+\ell/m}(X_0)\vee 1, 
\end{equation}
under the standing assumption   (\ref{aslambda}) of the Theorem.  Indeed, by virtue of (\ref{L}) and a little more precisely,  
\begin{eqnarray*}
& \displaystyle {\mathbb E}_{X, \mu} \bar\tau_{0011}^{k+1} = \int {\mathbb E}_{X, Y} \bar\tau_{0011}^{k+1} \mu(dY) 
\le C \int \bar V_{m,a+\ell/m}(X,Y) \vee 1\mu(dY) 
 \\\\
& \displaystyle = C V_{m,a+\ell/m}(X)\vee 1 + C \int V_{m,a+\ell/m}(Y)\vee 1\,\mu(dY) 
\le C (V_{m,a+\ell/m}(X)\vee 1 + C_1 + 1).
\end{eqnarray*}
Note that here both constants $C$ and $C_1$ admit  some effective estimates. 

~

\label{step10}
\noindent
{\bf \ref{step10}}. 
The final step of the proof is quite standard in coupling.
Consider two independent versions $X$ and $Y$ of
our Markov process, one
starting at $X_0$ and another at the stationary distribution~$\mu$ found earlier. 
Now, on some new probability space as described, for example, in \cite{But_Ver}, we estimate, uniformly in $A \in {\cal B}({\mathbb R}^1)$, 
\begin{eqnarray*}
& \displaystyle |(\mu^{X_0}_t - \mu)(A)| \le |{\mathbb E}_{X_0, \mu}(1(X_t\in A) -
1(Y_t\in A))| \, 1(t \ge \bar \tau_{0011})
 \\\\
& \displaystyle + |{\mathbb E}_{X_0, \mu}(1(X_t\in A) -
1(Y_t\in A))| \, 1(t < \bar \tau_{0011})
 \\\\
& \displaystyle \le {\mathbb E}_{X_0, \mu}1(t < \bar \tau_{0011}) ={\mathbb P}_{X_0, \mu}(t < \bar \tau_{0011})
\le \frac{{\mathbb E}_{X_0, \mu}\bar \tau_{0011}^{k+1}}{t^{k+1}}
\le \frac{C(V_{m,a+\ell/m}(X_0)\vee 1+\tilde C_1)}{t^{k+1}}
\end{eqnarray*}
(with \(\tilde C_1 = C_1+1\)). This shows the main statement  of the Theorem;
the version with the denominator \((1+t)^{k+1}\) clearly can be used instead of \(t^{k+1}\) since the left hand side may never exceed one. Uniqueness of a stationary distribution \(\mu\)   automatically follows from this  convergence. The Theorem \ref{thm1} is proved. 

\section*{Acknowledgements}
The author is sincerely grateful to  G.A. Zverkina for very useful discussions and for a fruitful criticism and to the anonymous referee for useful remarks and suggestions.

\end{document}